\documentclass{amsart}

\usepackage{indentfirst}
\usepackage{amssymb}
\usepackage{amsmath}
\usepackage{graphics}
\usepackage{epsfig}
\usepackage{MnSymbol,wasysym}

\usepackage[utf8]{inputenc}
\usepackage{latexsym}
\usepackage{graphicx}

\usepackage[normalem]{ulem}

\usepackage{hyperref}

\usepackage[T1]{fontenc}
\usepackage{kantlipsum}
\usepackage[usenames,dvipsnames]{xcolor}
\usepackage[breakable, theorems, skins]{tcolorbox}
\tcbset{enhanced}

\definecolor{light-gray}{gray}{0.95}

\catcode`@=11 \@addtoreset{equation}{section}
\renewcommand\theequation{\thesection.\@arabic\c@equation}
\catcode`@=12

\newcommand{\RR}{\mathbb{R}}

\expandafter\chardef\csname pre amssym.def
at\endcsname=\the\catcode`\@ \catcode`\@=11
\def\undefine#1{\let#1\undefined}
\def\newsymbol#1#2#3#4#5{\let\next@\relax
 \ifnum#2=\@ne\let\next@\msafam@\else
 \ifnum#2=\tw@\let\next@\msbfam@\fi\fi
 \mathchardef#1="#3\next@#4#5}
\def\mathhexbox@#1#2#3{\relax
 \ifmmode\mathpalette{}{\m@th\mathchar"#1#2#3}%
 \else\leavevmode\hbox{$\m@th\mathchar"#1#2#3$}\fi}
\def\hexnumber@#1{\ifcase#1 0\or 1\or 2\or 3\or 4\or 5\or 6\or 7\or 8\or
 9\or A\or B\or C\or D\or E\or F\fi}

\font\teneufm=eufm10 \font\seveneufm=eufm7 \font\fiveeufm=eufm5
\newfam\eufmfam
\textfont\eufmfam=\teneufm \scriptfont\eufmfam=\seveneufm
\scriptscriptfont\eufmfam=\fiveeufm

\catcode`\@=\csname pre amssym.def at\endcsname

\newcommand{\eqn}{\begin{eqnarray}}
\newcommand{\een}{\end{eqnarray}}

\newtheorem {Theorem}  {Theorem}

\numberwithin{Theorem}{section}
\newtheorem{Lemma}[Theorem]{Lemma}
\newtheorem{Proposition}[Theorem]{Proposition}
\newtheorem{Definition}[Theorem]{Definition}
\newtheorem{Remark}[Theorem]{Remark}

\newcommand{\NN}{{\mathbb N}}

\renewcommand{\a}{\alpha}

\renewcommand{\div}{\mbox{div}}

\begin{document}

\title[Decay of Navier-Stokes and Navier-Stokes-Coriolis in critical spaces]{Algebraic decay rates for 3D Navier-Stokes and Navier-Stokes-Coriolis equations in $ \dot{H}^{\frac{1}{2}}$}

\author[M. Ikeda]{Masahiro Ikeda}
\address[M. Ikeda]{Department of Mathematics, Faculty of Science and Technology, Keio University, 3-14-1 Hiyoshi, Kohoku-ku, Yokohama, 223-8522, Japan/Center for Advanced Intelligence Project, RIKEN, Japan}
\email{masahiro.ikeda@keio.jp/masahiro.ikeda@riken.jp}

\author[L. Kosloff]{Leonardo Kosloff}
\address[L. Kosloff]{Departamento de Matem\'atica, Instituto de Matem\'atica, Estat\'{\i}stica e
Computa\c{c}\~ao Cient\'{\i}fica, Universidade Estadual de Campinas, Rua Sergio Buarque de Holanda, 651, 13083-859, Campinas - SP, Brazil}
\email{kosloff@ime.unicamp.br}

\author[C. J. Niche]{C\'esar J. Niche}
\address[C.J. Niche]{Departamento de Matem\'atica Aplicada, Instituto de Matem\'atica. Universidade Federal do Rio de Janeiro, CEP 21941-909, Rio de Janeiro - RJ, Brazil}
\email{cniche@im.ufrj.br}

\author[G. Planas]{Gabriela Planas}
\address[G. Planas]{Departamento de Matem\'atica, Instituto de Matem\'atica, Estat\'{\i}stica e
Computa\c{c}\~ao Cient\'{\i}fica, Universidade Estadual de Campinas, Rua Sergio Buarque de Holanda, 651, 13083-859, Campinas - SP, Brazil}
\email{gplanas@unicamp.br}

\thanks{M. Ikeda has been supported by JSPS KAKENHI Grant Number JP 23K03174. L. Kosloff has been supported by FAPESP-Brazil  grant 2016/15985-0.  C.J. Niche acknowledges support from Bolsa PQ CNPq - 308279/2018-2  and PROEX - CAPES. G. Planas was partially supported by CNPq-Brazil grant 310274/2021-4, and FAPESP-Brazil grant 19/02512-5}



\date{\today}

\begin{abstract}
An algebraic upper bound for the decay rate of solutions to the Navier-Stokes and Navier-Stokes-Coriolis equations in the critical space $\dot{H} ^{\frac{1}{2}} (\RR ^3)$ is derived using the Fourier Splitting Method.  Estimates are framed in terms of the decay character of initial data, leading to solutions with algebraic decay and showing in detail the roles played by the linear and nonlinear parts.
\end{abstract}

\maketitle

\section{Introduction}
 \subsection{Critical spaces}  Our main goal in this article is to prove algebraic decay for solutions to the Navier-Stokes equations in $\RR ^3$
\begin{align}
\label{eqn:NS}
\partial _t u + (u \cdot \nabla) u +  \nabla p  & =  \Delta u , \nonumber \\ div \, u & =  0,  \\  u(x,0) & = u_0 (x), \nonumber
\end{align}
and for the Navier-Stokes-Coriolis equations in $\RR ^3$
\begin{align}
\label{eqn:navier-stokes-coriolis-omega}
\partial _t u + (u \cdot \nabla) u +  \nabla p + \Omega \, e_3 \times u & =  \Delta u , \nonumber \\ div \, u & =  0,  \\ u(x,0) & = u_0 (x), \nonumber
\end{align}
in the critical space $\dot{H} ^{\frac{1}{2}}(\RR ^3)$.  Solutions to  \eqref{eqn:NS} and \eqref{eqn:navier-stokes-coriolis-omega} have a natural scaling and translation symmetry
\begin{displaymath}
u _{\lambda, x_0} (x,t) = \lambda u (\lambda \left(x - x_0 \right), \lambda^2 t), \qquad \lambda > 0, x_0 \in \RR ^3
\end{displaymath}
i.e. if $u$ is a solution,  so is $u_{\lambda, x_0}$.  A Banach space $X$ is critical if $\Vert u_{\lambda, x_0 } \Vert _X = \Vert u \Vert _X$, for all $\lambda,  x_0 $.  As examples of critical  spaces in $\RR ^n$ for  \eqref{eqn:NS} and \eqref{eqn:navier-stokes-coriolis-omega} we mention
\begin{displaymath}
\dot{H} ^{\frac{n}{2} -1}; \, L^n; \, \dot{B} ^{-1 + \frac{n}{p}} _{p, \infty}, p \geq n; \, BMO ^{-1}; \dot{B} ^{-1} _{\infty, \infty},
\end{displaymath}
where $ \dot{B} ^{-1} _{\infty, \infty}$ is maximal,  in the sense that every $X$ that is critical  is continuously embedded in $ \dot{B} ^{-1} _{\infty, \infty}$.  An important feature of these spaces is that global mild solutions can be obtained for small enough initial data through Banach's Fixed Point Theorem and these solutions are well posed in the Hadamard sense,  with the notable exception of,  precisely,  $\dot{B} ^{-1} _{\infty, \infty}$, see Bourgain and Pavlovi\'{c} \cite{MR2473255},  Germain \cite{MR2473256}.   There is a vast amount of bibliography concerning results about solutions to \eqref{eqn:NS} and \eqref{eqn:navier-stokes-coriolis-omega} in critical spaces,  we refer to the books and surveys by Bahouri, Chemin and Danchin \cite{MR2768550}, Cannone \cite{MR2099035}, Gallagher \cite{MR3916784} and Lemari\'e-Rieusset \cite{MR1938147}, \cite{MR3469428} as a starting point.

\subsection{Decay results for Navier-Stokes equations} In order to provide context,  we briefly recall some  well known results on the decay of solutions to the Navier-Stokes equations \eqref{eqn:NS} in $L^2(\RR^n)$ and $H^s(\RR^n), s > 0$.    The formal equality
\begin{equation}
\label{eqn:energy-equality}
\frac{1}{2} \frac{d}{dt} \Vert u(t) \Vert _{L^2} ^2 = - \Vert \nabla u(t) \Vert _{L^2} ^2,
\end{equation}
led Leray \cite{MR1555394} to ask whether $\Vert u(t) \Vert _{L^2}$ goes to zero or not when time goes to infinity. Kato \cite{MR760047}  and Masuda \cite{MR767409} showed that indeed the $L^2$-norm goes to zero,  but they did not provide a rate of decay.  An algebraic decay rate was obtained by M.E. Schonbek \cite{MR775190}, \cite {MR837929}: for initial data $u_0 \in L^p (\RR ^n) \cap L^2 (\RR ^n)$, with $1 \leq p < 2$ she proved that
\begin{equation}
\label{eqn:lp-l2-decay}
\Vert u(t) \Vert _{L^2} ^2 \leq C (1 + t) ^{- \frac{n}{2} \left(\frac{2}{p} - 1 \right)}, \qquad t > 0.
\end{equation}
The key tool for attaining this result is the Fourier Splitting Method,  which amounts to obtaining a differential inequality for the $L^2$-norm,  in which the right hand side is  the average of the solution $u$ in a small,  time dependent shrinking ball around the origin in frequency space.  We carefully describe the Fourier Splitting Method in Section \ref{fourier-splitting-section}.   Note that in \eqref{eqn:lp-l2-decay} the rate is determined by the $L^p$ part of the initial datum $u_0$.  Building on this idea,  Wiegner \cite{MR881519} proved  decay under hypotheses on the decay of the  linear part and of the external forces.   Later, M.E. Schonbek and Wiegner \cite{MR1396285} proved decay for the  $H^k(\RR ^n)$ norm,   for $k \in \NN$.   Using a rather different method based on  Gevrey estimates,  Oliver and Titi \cite{MR1749867} were able to prove bounds for the radius of analiticity of solutions and estimates for decay for the $\dot{H}^s (\RR^n)$ norm  with $s > 0$,  provided information concerning the $L^2 (\RR^n)$ decay is known.
 Bae and Biswas \cite{MR3457535} proved a $L^p$-version of the decay result in \cite{MR1749867}, however,  estimates for decay hold only for the $\dot{H}^s_p (\RR^3)$-norm
for $s > \frac{3}{p}-1$. So, for $ p=2$, the critical case $ s =\frac{1}{2}$ is not achieved. We stress that results in \cite{MR1749867} and \cite{MR3457535}
assume initial datum in $L^2(\RR^n)$.  We only assume initial datum in $\dot{H}^\frac{1}{2}(\RR^3)$, hence the decay for $s=\frac{1}{2}$ does not come out from \cite{MR1749867} in this case.

A more general approach was introduced by Bjorland and M.E. Schonbek \cite{MR2493562},  and later refined by Niche and M.E. Schonbek \cite{MR3355116}  and Brandolese \cite{MR3493117},  by associating to a given $u_0 \in L^2(\RR ^3)$ a decay character $r^{\ast} = r^{\ast} (u_0)$,  $- \frac{n}{2} < r^{\ast} < \infty$,  which roughly speaking says that $|\widehat{u_0} (\xi)| \approx |\xi| ^{r^{\ast}}$ when $|\xi| \approx 0$.  See Section \ref{DecayCharacter} for details on the decay character. Using this quantity,  it is possible to prove that for weak solutions to \eqref{eqn:NS}
\begin{equation}
\label{eqn:decay-with-decay-character}
\Vert u(t) \Vert _{L^2} ^2 \leq C (1 + t) ^{- \min \left\{ \frac{3}{2} + r^{\ast}, \frac{5}{2} \right\}}.
\end{equation}
From \eqref{eqn:decay-with-decay-character} we conclude that for initial data with $r^{\ast} \leq 1$,  the decay is driven by the slower linear part,  while for $r^{\ast} > 1$ the nonlinear part becomes slower.  See  Section \ref{DecayCharacter} for estimates on dissipative linear equations,  including the heat equation.

For a thorough survey on results concerning decay and long time behaviour of the Navier-Stokes equations,  see  Brandolese and M.E. Schonbek \cite{MR3916783}.

\subsection{Main results: algebraic decay in critical space $\dot{H} ^{\frac{1}{2}}(\RR ^3)$}
Unlike the case of the Navier-Stokes equations discussed in the previous Section, there are few
results on the decay of solutions to \eqref{eqn:NS} and   \eqref{eqn:navier-stokes-coriolis-omega} in the critical space
 $\dot{H} ^{\frac{1}{2}}(\RR ^3)$
 assuming that the  initial datum belongs just to the critical space $\dot{H}  ^{\frac{1}{2}} (\RR ^3)$. It  is our main goal in this article to prove algebraic decay  in  $\dot{H}  ^{\frac{1}{2}} (\RR ^3)$  for mild solutions to \eqref{eqn:NS} and   \eqref{eqn:navier-stokes-coriolis-omega} through the combined use of the Fourier Splitting Method and the decay character of initial data.

\subsubsection{Navier-Stokes equations} Well posedness of mild solutions to \eqref{eqn:NS} in $\dot{H} ^{\frac{1}{2}} (\RR ^3)$ is a classical result,  which can be traced back to Fujita and Kato \cite{MR166499} (see also Bahouri, Chemin and Danchin \cite{MR2768550},  and Lemari\'e-Rieusset \cite{MR1938147}, \cite{MR3469428}  for a textbook version).   Known decay estimates for solutions to \eqref{eqn:NS} merely provide decay to zero with no explicit rates,  i.e.
\begin{equation}
\label{eqn:decay-gallagher-iftimie-planchon}
\lim _{t \to \infty} \Vert u(t) \Vert _{\dot{H} ^{\frac{1}{2}}} = 0,
\end{equation}
as proved by Gallagher, Iftimie and Planchon \cite{MR1891005}.  Sawada \cite{MR2175199} proved space analyticity of mild solutions in $\dot{H} ^{\frac{n}{2} - 1} (\RR ^n)$ but his estimates provide decay for $\dot{H} ^{\frac{n}{2} - 1 + \alpha}$-norms only,  with $\alpha > 0$.  As $\dot{H} ^{\frac{1}{2}}(\RR ^n)$ is an $L^2$-based Sobolev space,  we expect the Fourier Splitting Method may allow us to obtain decay rates in \eqref{eqn:decay-gallagher-iftimie-planchon}.

We state now our first result.   From now on we denote by $\Lambda = (- \Delta) ^{\frac{1}{2}}$ the Riesz potential.

\begin{Theorem} \label{decay-ns-critical}   Let $u_0 \in \dot{H} ^{\frac{1}{2}} (\RR ^3)$, with $ div \, u_0 = 0 $ and $\Vert u_0  \Vert _{ \dot{H} ^{\frac{1}{2}}} < \epsilon$, for small enough $\epsilon >  0$ and  $-\frac{3}{2} < q^{\ast} = r^{\ast} (\Lambda ^{\frac{1}{2}} u_0) < \infty$. Then, for any  mild solution to \eqref{eqn:NS}, we have that for some constant $C=C(\epsilon,q^{\ast})$
 \begin{displaymath}
\Vert u (t) \Vert _{\dot{H} ^{\frac{1}{2}}} ^2  \leq C (1+t) ^{- \min \left\{ \frac{3}{2} + q^{\ast}, 1 \right\}}.
\end{displaymath}
\end{Theorem}

The starting point for this result is the rigorous proof of the inequality
\begin{equation}
\label{eqn:formal-inequality}
\frac{1}{2} \frac{d}{dt} \Vert u(t) \Vert^2 _{\dot{H} ^{\frac{1}{2}}} \leq - \left(1 -  C \Vert  u(t) \Vert _{\dot{H} ^{\frac{1}{2}}} \right) \Vert \nabla u(t) \Vert^2 _{\dot{H} ^{\frac{1}{2}}},
\end{equation}
which,  as $\Vert u(t) \Vert _{\dot{H} ^{\frac{1}{2}}}$ is small for all time due to the fact that solutions are obtained by a fixed point argument, leads to an inequality analogous to \eqref{eqn:energy-equality} but for the $\dot{H} ^{\frac{1}{2}}$-norm.  We then adapt the Fourier Splitting Method to this context and prove decay.  Note that, as the $L^2$ decay in \eqref{eqn:decay-with-decay-character}, the decay  in Theorem \ref{decay-ns-critical} is driven by  the linear and nonlinear parts for different sets of initial data.

Recently,  Duong,  Khai and Tri \cite{MR4052578} proved that for initial data in $L^3 (\RR ^3)$ such that
\begin{displaymath}
\Vert e^{t \Delta} u_0 \Vert _{L^3} = o(t^{-\alpha}),  \, \Vert e^{t \Delta} u_0 \Vert _{L^3} = O(t^{-\alpha}), \qquad 0 \leq \alpha \leq 1,
\end{displaymath}
then mild solutions to the Navier-Stokes equations are such that, respectively
\begin{displaymath}
\Vert u (t) \Vert _{L^3} = o(t^{-\alpha}),  \, \Vert u (t) \Vert _{L^3} = O(t^{-\alpha}), \qquad 0 \leq \alpha \leq 1.
\end{displaymath}
As $\dot{H} ^{\frac{1}{2}} (\RR ^3) \subset L^3 (\RR ^3)$, we have
\[ \Vert u (t) \Vert _{L^3}^2 \leq C (1+t) ^{- \min \left\{ \frac{3}{2} + q^{\ast}, 1 \right\}},
\]
for any  $ u_0 \in \dot{H} ^{\frac{1}{2}} (\RR ^3)$. Thus, this refines and extends their results.

\subsubsection{Navier-Stokes-Coriolis equations} In this setting, we are also motivated to study  the Navier-Stokes-Coriolis equations \eqref{eqn:navier-stokes-coriolis-omega}.
These equations are a prototype for geophysical models with strong rotation around a fixed axis, this being modelled by the last term in the left hand side,  the Coriolis term, which is typically larger than the other terms because of the parameter $\Omega \in \RR$.   For a thorough study of the mathematical theory of this equation and others modelling rotating fluids  see Chemin \emph{et al.} \cite{MR2228849}.  Strong rotation is essential for proving  global well-posedness in the case of large data in $H^{\frac{1}{2}} (\RR ^3)$, and in this framework it can also be shown that the asymptotic limit system is determined by the 2D Navier-Stokes equations, (see   Theorem 5.7 in Chemin {\em et al.}  \cite{MR2228849} for a precise statement).  For similar  results in the case of domains,  including the treatment of the Ekman boundary layer problem,  we refer to Chemin {\em et al.} \cite{MR2228849},  Giga,  Inui,  Mahalov and Matsui \cite{MR2342201}, and  Iftimie,  Raugel and Sell \cite{MR2333468}.

M.E. Schonbek and Vallis \cite{MR1689401} proved decay results for the $L^2$-norm of solutions to a system of coupled Boussinesq and generalized  damped Navier-Stokes-Coriolis equations.  However,  the proof of the $L^2$ decay of the velocity $u$  for this system relies on the damping and coupling to the temperature $\theta$.  As a result of this,  these results do not directly apply to  \eqref{eqn:navier-stokes-coriolis-omega}.  For mild solutions to the Navier-Stokes-Coriolis equations \eqref{eqn:navier-stokes-coriolis-omega} in $H ^{\frac{1}{2}} (\RR ^3)$  and in $\dot{H}  ^{\frac{1}{2}} (\RR ^3)$ only boundedness is known,  see Hieber and Shibata \cite{MR2609321} and Iwabuchi and Takada \cite{MR3096523}, respectively.  Recently, Ahn, Kim and Lee \cite{MR4444120} addressed decay of solutions to a fractional version of \eqref{eqn:navier-stokes-coriolis-omega}, but their results do not cover the critical case $\dot{H}^{\frac{1}{2}}(\RR^3)$. Egashira and Takada \cite{MR4544280} proved some decay rates for solutions to  \eqref{eqn:navier-stokes-coriolis-omega} in   $L^p$-norm, for $p \geq 2$, by considering initial datum in  $ \dot{H}^s(\RR^3) \cap L^1(\RR^3)$, with $s > \frac{1}{2}$, again not covering the critical case $s = \frac{1}{2}$.

We now state our result for  decay in the critical space $\dot{H}^{\frac{1}{2}} (\RR ^3)$ for  \eqref{eqn:navier-stokes-coriolis-omega}.

\begin{Theorem}  \label{thm-decay-h-dot-one-half} Let $u_0 \in \dot{H} ^{\frac{1}{2}} (\RR ^3)$, with $ div \, u_0 = 0 $ and $\Vert u_0  \Vert _{ \dot{H} ^{\frac{1}{2}} (\RR ^3)} < \epsilon$, for small enough $\epsilon > 0 $. Then, there exists $\omega(u_0) > 0$ such that for any $\Omega$ with $|\Omega| > \omega(u_0)$,  if $-\frac{3}{2} < q^{\ast} = r^{\ast} (\Lambda ^{\frac{1}{2}} u_0) < \infty$, then, for any mild  solution to \eqref{eqn:navier-stokes-coriolis-omega}, we have that for a constant $C=C(\epsilon,q^{\ast},|\Omega|)$
\begin{displaymath}
\Vert u (t) \Vert _{\dot{H} ^{\frac{1}{2}}} ^2  \leq C (1+t) ^{- \min \left\{ \frac{3}{2} + q^{\ast}, 1 \right\}}.
\end{displaymath}
\end{Theorem}

As in Theorem \ref{decay-ns-critical}, we turn \eqref{eqn:formal-inequality} into  an inequality analogous to \eqref{eqn:energy-equality} but for the
$\dot{H}^{\frac{1}{2}}$-norm,  which we do through the use of the smallness of the solution.  We also adapt to this context the Fourier Splitting Method.

\subsection{Organization of this article} In  Section \ref{preliminary}, we present some preliminary results: we review the Fourier Splitting Method and the concept of decay character, and we recall some important inequalities that  will be used in our proofs.
In Section \ref{NS} we  proceed to set up the semigroup framework in which the fixed point arguments yield the global existence of solutions to Navier-Stokes equations \eqref{eqn:NS} in the critical space $\dot{H}^{\frac{1}{2}}(\RR^3)$.  We conclude this Section with the proof of the decay estimate of Theorem \ref{decay-ns-critical}.
Finally, in Section \ref{NSC} we carry out the analysis of the Navier-Stokes-Coriolis equations leading to the proof of Theorem \ref{thm-decay-h-dot-one-half}.
 An additional decay estimate for Navier-Stokes-Coriolis equations \eqref{eqn:navier-stokes-coriolis-omega} in $H^{\frac{1}{2}}(\RR^3)$ is also provided.
The proof of this estimate is rather simple because in this case we do have $L^2$ decay estimates at hand, due to the fact that the initial datum is in $H^{\frac{1}{2}} (\RR^3) \subset L^2 (\RR^3)$.

\subsection{Acknowledgments} We thank Edriss Titi for comments and remarks that improved the presentation of this article and for pointing out references related to our work.

\section{Preliminary results}
\label{preliminary}

\subsection {Fourier Splitting} \label{fourier-splitting-section}The Fourier Splitting Method was developed by M.E. Schonbek  to study decay of energy for solutions to parabolic conservations laws \cite{MR571048} and to Navier-Stokes equations \cite{MR775190}, \cite{MR837929} and rests on the  observation that for these equations ``for large enough times, remaining energy is concentrated at the low frequencies''. We describe now a rather general setting, introduced by Niche and M.E. Schonbek \cite{MR3355116}, in which this method can be used. Let $X$ be a Hilbert space and consider a linear diagonalizable pseudodifferential operator $\mathcal{L}: X^n \to \left( L^2 (\RR^n) \right) ^n$, this is an $\mathcal{L}$ which has an associated symbol $ M(\xi)$ such that
\begin{equation}
\label{eqn:symbol}
M(\xi) = P^{-1} (\xi) D(\xi) P(\xi), \qquad \xi-a.e.
\end{equation}
where $P(\xi) \in O(n)$ and $D(\xi) = - c_i |\xi|^{2\a} \delta _{ij}$, for $c_i > c>0$ and $0 < \a \leq 1$. Take the linear system
\begin{equation}
\label{eqn:linear-part}
v_t = \mathcal{L} v,
\end{equation}
multiply by $v$ and obtain
\begin{align}
\frac{1}{2} \frac{d}{dt} \Vert v(t) \Vert _{L^2} ^2 & =   \langle \widehat{v}, M \widehat{ v} \rangle _{L^2} =  \langle \widehat{v}, P^{-1}  D P \widehat{v} \rangle _{L^2} \nonumber \\ & =  - \langle (-D) ^{\frac{1}{2}} P \widehat{v} ,  (-D) ^{\frac{1}{2}} P \widehat{v} \rangle _{L^2} \nonumber \\ & =  - \int _{\RR^n} |(-D) ^{\frac{1}{2}} P \widehat{v}|^2 \, d \xi \nonumber \\ & \leq  - C  \int _{\RR^n} |\xi|^{2 \a} |\widehat{v}|^2 \, d \xi. \notag
\end{align}
Thus, the $L^2$-norm decreases.

\begin{Remark}
We note that the vectorial fractional Laplacian $(- \Delta) ^{\alpha}$, with $0 < \alpha \leq 1$; the vectorial Lam\'e  operator
$\mathcal{L} u = \Delta u + \nabla \, \div \, u, $
see  Example 2.9 in Niche and M.E. Schonbek \cite{MR3355116}, and the linear part of the magneto-micropolar system (see Niche and Perusato \cite{MR4379088}), amongst others, are diagonalizable as in \eqref{eqn:symbol}.
\end{Remark}

Now, we split frequency space into two parts by introducing a ball around the origin $B(t)$, with decreasing time dependent radius. More precisely, take $B(t) = \{\xi \in \RR^n: |\xi| \leq g(t) \}$, for a nonincreasing, continuous $g$ to be determined later, with $g(0) >0$. Then, for $\alpha =1$
\begin{align*}
\frac{d}{dt} \Vert \widehat{v} (t) \Vert _{L^2} ^2 & \leq - 2C  \int _{\RR ^n} |\xi|^2 |\widehat{v} (\xi,t)|^2 \, d \xi  \leq -2C \int _{B(t) ^{c}} |\xi|^2 |\widehat{v} (\xi,t)|^2 \, d \xi \notag  \\ & \leq - 2C g^2 (t)  \int _{B(t) ^{c}}  |\widehat{v} (\xi,t)|^2 \, d \xi \\ & = - 2C g^2 (t) \Vert \widehat{v} (t)  \Vert _{L^2} ^2 \notag + 2C g^2 (t)  \int _{B(t)}  |\widehat{v} (\xi,t)|^2 \, d \xi, \notag
\end{align*}
which, after using the integrating factor
\begin{displaymath}
h(t) = \exp \left( \int _0 ^t 2C  g^2 (s) \, ds \right),
\end{displaymath}
leads to the key inequality
\begin{equation}
\label{eqn:key-inequality-clean}
\frac{d}{dt}   \left( h(t)  \Vert v (t) \Vert _{L^2} ^2 \right) \leq h'(t)  \int _{B(t)} |\widehat{v} (\xi, t)| ^2 \, d \xi.
\end{equation}
Now,  we would have to appropriately choose $g$ and find a pointwise estimate for $|\widehat{v} (\xi, t)|$ to obtain an upper bound after integrating \eqref{eqn:key-inequality-clean}.

Note that because of the form of our linear operator, i.e. \eqref{eqn:symbol},  the linear part of $\widehat{v}$ will be such that
\begin{equation}
\label{eqn:large-time}
e^{t M(\xi)} \, \widehat{v}_0 (\xi) \approx \widehat{v}_0 (\xi), \qquad t \gg 1,
\end{equation}
hence the  decay of the linear part will be essentially determined by $\widehat{v}_0 (\xi)$, with $\xi \approx 0$. This motivates the definition of decay character in the next Section.

\subsection{Decay character} \label{DecayCharacter} We now recall the definition and properties of the decay character, introduced by Bjorland and M.E. Schonbek \cite{MR2493562}, and refined by Niche and M.E. Schonbek \cite{MR3355116}, and Brandolese \cite{MR3493117}. Loosely speaking, to any initial datum  $v_0 \in L^2(\RR^n)$ the decay character associates a number which measures its ``algebraic order'' near the origin, comparing $ |\widehat{v_0} (\xi)|$ to $f(\xi) = |\xi|^{r}$ at $\xi = 0$. Due to \eqref{eqn:large-time}, $|\widehat{v_0} (\xi)|^2$  determines the behaviour of the right hand side of \eqref{eqn:key-inequality-clean}, thus providing decay estimates.

\begin{Definition} \label{decay-indicator}
Let  $v_0 \in L^2(\RR^n)$. For $r \in \left(- \frac{n}{2}, \infty \right)$, we define the {\em decay indicator}  $P_r (v_0)$ corresponding to $v_0$ as
\begin{displaymath}
P_r(v_0) = \lim _{\rho \to 0} \rho ^{-2r-n} \int _{B(\rho)} \bigl |\widehat{v_0} (\xi) \bigr|^2 \, d \xi,
\end{displaymath}
provided this limit exists. In the expression above,  $B(\rho)$ denotes the ball at the origin with radius $\rho$.
\end{Definition}

\begin{Definition} \label{df-decay-character} The {\em decay character of $ v_0$}, denoted by $r^{\ast} = r^{\ast}( v_0)$ is the unique  $r \in \left( -\frac{n}{2}, \infty \right)$ such that $0 < P_r (v_0) < \infty$, provided that this number exists. We set $r^{\ast} = - \frac{n}{2}$, when $P_r (v_0)  = \infty$ for all $r \in \left( - \frac{n}{2}, \infty \right)$  or $r^{\ast} = \infty$, if $P_r (v_0)  = 0$ for all $r \in \left( -\frac{n}{2}, \infty \right)$.
\end{Definition}

It is possible to explicitly compute the decay character for many important examples. When $v_0 \in L^p (\RR^n) \cap L^2 (\RR ^n)$ for $1 < p < 2$ and $v_0 \notin L^{\bar{p}} (\RR ^n)$ for $\bar{p} < p$, we have that $r^{\ast} (v_0) = - n \left( 1 - \frac{1}{p} \right)$, see Example 2.6 in Ferreira, Niche and Planas \cite{MR3565380}. When $v_0 \in L^1 (\RR^n) \cap L^2 (\RR^n)$ and $|\widehat{v_0} (\xi)|$ is bounded away from zero near the origin, i.e.
\begin{displaymath}
0< C_1 \leq |\widehat{v_0} (\xi)| \leq C_2, \qquad |\xi| \leq \beta,
\end{displaymath}
for some $\beta > 0$ and $0 < C_1 \leq C_2$, then $r^{\ast} (v_0) = 0$, see Section 4 in M.E. Schonbek \cite{MR837929}. For  $v_0 \in L^{1, \gamma} (\RR^n)$, for $0 \leq \gamma \leq 1$, i.e.
\begin{equation*}
%\label{eqn:l-uno-gamma}
\Vert f \Vert _{L^{1, \gamma} (\mathbb{R} ^n)} = \int _{\mathbb{R} ^n} \left( 1 + |x| \right) ^{\gamma} |f(x)| \, dx < \infty,
\end{equation*}
we have from Lemma 3.1 in Ikehata \cite{MR2055280} that

\begin{displaymath}
|\widehat{v_0} (\xi) | \leq C |\xi| ^{\gamma} + \left| \int _{\RR ^n} v_0(x) \, dx \right|, \qquad \xi \in \RR ^n.
\end{displaymath}
If $v_0$ has zero mean,  then $r^{\ast} (v_0) = \gamma$, if not, then $r^{\ast} (v_0) = 0$.

We now state the Theorem that describes decay in terms of the decay character for linear operators as in \eqref{eqn:symbol}.

\begin{Theorem}{(Theorem 2.10, Niche and M.E. Schonbek \cite{MR3355116})}
\label{characterization-decay-l2}
Let $v_0 \in L^2 (\RR^n)$ have decay character $r^{\ast} (v_0) = r^{\ast}$. Let $v (t)$ be a solution to  \eqref{eqn:linear-part} with initial datum $v_0$, where the operator $\mathcal{L}$ is such that \eqref{eqn:symbol}  holds. Then if $- \frac{n}{2 } < r^{\ast}< \infty$, there exist constants $C_1, C_2> 0$ such that
\begin{displaymath}
C_1 (1 + t)^{- \frac{1}{\a} \left( \frac{n}{2} + r^{\ast} \right)} \leq \Vert v(t) \Vert _{L^2} ^2 \leq C_2 (1 + t)^{- \frac{1}{\a} \left( \frac{n}{2} + r^{\ast} \right)}.
\end{displaymath}
\end{Theorem}

Let $\Lambda = (- \Delta) ^{\frac{1}{2}}$ and $u_0 \in H^s(\RR^n)$, for $s > 0$.  We will  need the following relation between the decay character $r^{\ast} _s (u_0) = r^{\ast} \left( \Lambda^s u_0 \right)$ and that of $u_0$, $r^{\ast} (u_0)$.

\begin{Theorem}{(Theorem 2.11, Niche and M.E. Schonbek \cite{MR3355116})}
%\label{decay-character-hs}
Let $u_0 \in H^s (\RR^n)$, $s > 0$. Then, if    $-\frac{n}{2} < r^{\ast} (u_0) < \infty$ we have that
\begin{displaymath}
- \frac{n}{2} +s< r_s^{\ast}(u_0) < \infty, \quad r_s^{\ast}(u_0) = s + r^{\ast} (u_0).
\end{displaymath}
\end{Theorem}

As a consequence of this, we obtain the following decay estimates for homogeneous Sobolev norms.
\begin{Theorem}{(Theorem 2.12, Niche and M.E. Schonbek \cite{MR3355116})}
\label{characterization-decay-hs}
Let $v_0 \in H^s(\RR^n)$, $s > 0$ have decay character $r^{\ast} _s = r^{\ast} _s (v_0)$. Then if $- \frac{n}{2 } < r^{\ast} < \infty$, there exist constants $C_1, C_2 > 0$ such that
\begin{displaymath}
C_1 (1 + t)^{- \frac{1}{\a} \left( \frac{n}{2} + r^{\ast}  + s \right)} \leq \Vert v(t) \Vert _{\dot{H}^s} ^2 \leq C_2 (1 + t)^{- \frac{1}{\a} \left( \frac{n}{2} + r^{\ast}  + s \right)}.
\end{displaymath}
\end{Theorem}

 In Definitions \ref{decay-indicator} and \ref{df-decay-character}, it is assumed that those limits exist, thus giving rise to a positive $P_r (u_0)$. However this does not hold for all of $v_0 \in L^2 (\RR ^n)$, as Brandolese \cite{MR3493117} constructed initial data, highly oscillating near the origin, for which the limit in Definition \ref{decay-indicator} does not exist for some $r$. As a result of this, the decay character does not exist.  Brandolese also gave a slightly different definition of decay character, more general than that in  Definitions \ref{decay-indicator} and \ref{df-decay-character}, but which produces the same result when these hold. He also proved that the decay character $r ^{\ast}$ (in his more general version) exists for $v_0 \in L^2 (\RR ^n)$ if and only if $v_0$ belongs to a specific subset of a homogeneous Besov space, i.e. $v_0  \in  \dot{\mathcal{A}} ^{- \left(\frac{n}{2} + r^{\ast} \right)} _{2, \infty} \subset \dot{B} ^{- \left(\frac{n}{2} + r^{\ast} \right)} _{2, \infty}$. Moreover, for diagonalizable linear operators $\mathcal{L}$ as in \eqref{eqn:symbol}, solutions to the linear system \eqref{eqn:linear-part} with initial data $v_0$ have algebraic decay, i.e.
\begin{displaymath}
C_1 (1 + t)^{- \frac{1}{\a} \left( \frac{n}{2} + r^{\ast} \right)} \leq \Vert v(t) \Vert _{L^2} ^2 \leq C_2 (1 + t)^{- \frac{1}{\a} \left( \frac{n}{2} + r^{\ast} \right)},
\end{displaymath}
if and only if the decay character $r^{\ast} = r^{\ast} (v_0)$ exists. This provides a complete and sharp characterization of algebraic decay rates for such systems and provides a key tool for studying decay for nonlinear systems.

\subsection{Inequalities} We will need the following two versions of Gronwall's inequality.

\begin{Proposition}[Theorem 1, page 356, Mitrinovi\'{c}, Pe\v{c}ari\'{c} and Fink \cite{MR1190927}] \label{gronwall-1} Let $x,k:J \to \RR$ continuous and $a,b: J \to \RR$ Riemann integrable in $J = [\alpha, \beta]$. Suppose that $b, k \geq 0$ in $J$. Then, if
\begin{displaymath}
x(t) \leq a(t) + b(t) \int _{\alpha} ^t k(s) x(s) \, ds, \quad t \in J
\end{displaymath}
then
\begin{displaymath}
x(t) \leq a(t) + b(t) \int _{\alpha} ^t a(s) k(s) \exp \left( \int_s ^t b(r)k(r)  \, dr \right) \, ds, \quad t \in J.
\end{displaymath}
\end{Proposition}

\begin{Proposition}[Corollary 1.2,  page 4,  Bainov and Simeonov \cite{MR1171448}] \label{gronwall-2} Let $a,k,u,:J \to \RR$ continuous in $J = [\alpha, \beta]$ and $ k \geq 0$.  If $a(t)$ is nondecreasing then
\begin{displaymath}
\psi(t) \leq a(t) +  \int _{\alpha} ^t k(s) \psi(s) \, ds, \quad t \in J
\end{displaymath}
implies
\begin{displaymath}
\psi(t) \leq a(t)  \exp \left( \int_{\alpha} ^t k(s)  \, ds \right) \quad t \in J.
\end{displaymath}
\end{Proposition}

\section{Navier-Stokes equations}
\label{NS}

 We start recalling  results concerning solutions to Navier-Stokes equations \eqref{eqn:NS} that we will need in this work. We follow Chapter 5 of Bahouri {\em et al.}  \cite{MR2768550}. For other important sources of information for the incompressible Navier-Stokes equations in critical spaces, see Cannone \cite{MR2099035},  Gallagher \cite{MR3916784},  and Lemari\'e-Rieusset \cite{MR1938147}, \cite{MR3469428}.

Given a vector field $v$ in $\RR^3$, through the Leray projector $\mathbb{P}$ we obtain a divergence-free vector field $\mathbb{P} v$ with components
\begin{equation*}
\widehat{\left( \mathbb{P} v \right)} ^j (\xi) = \sum _{k = 1} ^3 \left(  \delta_{jk} -   \frac{\xi _j \xi_k}{|\xi|^2} \right)  \widehat{v^k} (\xi).
\end{equation*}

A mild solution, or simply solution, to \eqref{eqn:NS} is a function defined on a Banach space $E$, which is a fixed point of the map $\mathcal{G}: E \to E$ given by
\begin{equation}
\label{eqn:fixed-point-operator}
\mathcal{G} (u) = e^{t \Delta} u_0 -  \int_0 ^t e^{(t - s) \Delta} \mathbb{P} \,  \nabla \cdot \left( u \otimes u \right) (s) \, ds.
\end{equation}
From Theorem 5.6 in Bahouri {\em et al.}  \cite{MR2768550}, we know that there exists a (small) $\epsilon > 0$ such that if $\Vert u_0 \Vert _{\dot{H} ^{\frac{1}{2}}} < \epsilon$ then the $3D$ Navier-Stokes equations \eqref{eqn:NS} have a unique solution in $E = L^4 \bigl( [0, \infty); \dot{H} ^1  (\RR ^3) \bigr)$,  which is also in $F = C\bigl( [0, \infty); \dot{H} ^{\frac{1}{2}} (\RR ^3)  \bigr) \cap L^2 \bigl( [0, \infty); \dot{H} ^{\frac{3}{2}} (\RR ^3) \bigr)$. As the solution is obtained through Banach's Fixed Point Theorem, we have that the smallness condition $\Vert u_0 \Vert _{\dot{H} ^{\frac{1}{2}}} < \epsilon$ leads to
\begin{equation}
\label{eqn:small-from-start}
\Vert u (t) \Vert _{\dot{H} ^{\frac{1}{2}}} < 2 \epsilon, \qquad \forall \, t >0.
\end{equation}

We now recall the following Lemma.
\begin{Lemma}[Lemma 5.10, Bahouri {\em et al.}  \cite{MR2768550}]  Let $v \in C\left([0,T]; \mathcal{S} '(\RR ^n) \right)$ be the solution to
\begin{align*}
\partial _t v - \Delta v  & =   f,  \nonumber \\   v(x,0) & = v_0 (x),
\end{align*}
where $f \in L^2 \left([0,T]; \dot{H} ^{s-1}(\RR^n) \right)$ and $v_0 \in \dot{H}^{s} (\RR^n)$.
Then we have the energy identity
\begin{equation*}
%\label{eqn:energy-identity}
\Vert v(t) \Vert ^2 _{\dot{H} ^s}  + 2 \int _0 ^t \Vert \nabla v(\tau) \Vert ^2  _{\dot{H}^s} \, d\tau = \Vert v_0 \Vert ^2 _{\dot{H} ^s} + 2  \int _0 ^t \langle v (\tau), f(\tau) \rangle _{\dot{H}^s} \, d \tau.
\end{equation*}
\end{Lemma}

As an immediate consequence of this Lemma, solutions to \eqref{eqn:NS} obey the energy identity  with $f = - \mathbb{P} \nabla\cdot\left ( u \otimes u \right)(\tau)$
\begin{align*}
%\label{eqn:energy-identity-navier-stokes}
\Vert u(t) \Vert ^2 _{\dot{H} ^{\frac{1}{2}}} + 2 \int _0 ^t \Vert \nabla u(\tau) \Vert ^2  _{\dot{H}^{\frac{1}{2}}} \, d\tau &= \Vert v_0 \Vert ^2 _{\dot{H} ^{\frac{1}{2}}} \notag \\ & - 2  \int _0 ^t \langle  u(\tau), \mathbb{P} \nabla\cdot \left ( u \otimes u \right)(\tau) \rangle _{\dot{H}^{\frac{1}{2}}} \, d \tau.
\end{align*}
This identity is the key  for the proof, in Proposition 5.13 in Bahouri {\em et al.}  \cite{MR2768550}, that the $\dot{H} ^{\frac{1}{2}}$ norm is a Lyapunov function, i.e. it is a decreasing function of time on non-stationary solutions. Moreover, as we mentioned in the Introduction,
\begin{displaymath}
\lim _{t \to \infty} \Vert u (t) \Vert _{\dot{H} ^{\frac{1}{2}}} = 0,
\end{displaymath}
see Gallagher {\em et al.}  \cite{MR1891005} (and also Theorem 5.17 in Bahouri {\em et al.}  \cite{MR2768550}).

In order to obtain energy inequalities for equation \eqref{eqn:NS}  we will need the  estimate
\begin{align}
\label{eqn:bilinear}
\langle  u ,  \mathbb{P} \, \nabla \cdot ( u \otimes u)    \rangle _{\dot{H} ^{\frac{1}{2}}} & = \langle \Lambda ^{\frac{1}{2}} u  , \Lambda ^{\frac{1}{2}} \mathbb{P} \, \nabla \cdot ( u \otimes u)    \rangle _{L^2}  \notag \\ & \leq   C \Vert  u \Vert^2_{\dot{H} ^1}  \Vert \nabla u \Vert_{\dot{H} ^{\frac{1}{2}}}  \leq C \Vert  u \Vert_{\dot{H} ^{\frac{1}{2}}}  \Vert \nabla u \Vert^2 _{\dot{H} ^{\frac{1}{2}}} ,
\end{align}
where we used Lemma 5.12 in Bahouri \emph{et al.}  \cite{MR2768550} and interpolation.

\subsection{Proof of Theorem \ref{decay-ns-critical}.} We know that the $\dot{H} ^{\frac{1}{2}}$-norm is a nonincreasing function, hence has a derivative a.e. Then
\begin{align}
\notag \frac{d}{dt} \Vert u(t) \Vert^2 _{\dot{H} ^{\frac{1}{2}}} & = 2 \langle \Lambda ^{\frac{1}{2}} u(t), \partial_t \Lambda ^{\frac{1}{2}} u(t) \rangle \\ \notag & = 2 \langle \Lambda ^{\frac{1}{2}} u (t) ,\Lambda ^{\frac{1}{2}} \left( \Delta u (t)  - ( \mathbb{P} \, \nabla \cdot ( u \otimes u)  (t) \right)   \rangle \\ \notag &
 = -2 \Vert \nabla u (t) \Vert^2 _{\dot{H} ^{\frac{1}{2}}} - 2 \langle \Lambda ^{\frac{1}{2}} u (t) , \Lambda ^{\frac{1}{2}} \mathbb{P} \, \nabla \cdot ( u \otimes u)  (t)  \rangle. \notag
\end{align}
Now, after using \eqref{eqn:bilinear} in the last term there follows
\begin{displaymath}
\frac{1}{2} \frac{d}{dt} \Vert u(t) \Vert^2 _{\dot{H} ^{\frac{1}{2}}} \leq - \left(1 -  C \Vert  u(t) \Vert _{\dot{H} ^{\frac{1}{2}}} \right) \Vert \nabla u(t) \Vert^2 _{\dot{H} ^{\frac{1}{2}}},
\end{displaymath}
which is a well known estimate,  already formally obtained by  Kato \cite{MR1084601} and Cannone, Section 7.1 in \cite{MR2099035}.  From \eqref{eqn:small-from-start}, we know that  $\Vert u_0 \Vert _{\dot{H} ^{\frac{1}{2}}} < \epsilon$ leads to $\Vert u (t) \Vert _{\dot{H} ^{\frac{1}{2}}} < 2 \epsilon$, for all $t > 0$, hence choosing $\epsilon < \frac{1}{C}$ we have
\begin{displaymath}
\frac{d}{dt} \Vert u(t) \Vert^2 _{\dot{H} ^{\frac{1}{2}}} \leq - 2  \Vert \nabla u(t) \Vert^2 _{\dot{H} ^{\frac{1}{2}}},
\end{displaymath}
and then we can use the Fourier Splitting Method from Section \ref{fourier-splitting-section}. We thus obtain, for a nonincreasing and continuous $g$ to be determined later, with $g(0) = K > 0$, the inequality
\begin{align}
\label{eqn:key-inequality-ns}
\frac{d}{dt} & \left( \exp \left( \int _0 ^t 2  g^2 (s) \, ds \right)  \Vert \Lambda ^ {\frac{1}{2}} u (t) \Vert _{L^2} ^2 \right)   \notag \\ & \leq 2 g^2 (t) \left( \exp \left( \int _0 ^t 2  g^2 (s) \, ds \right) \right) \int _{B(t)} ||\xi| ^{\frac{1}{2}} \widehat{u} (\xi, t)| ^2 \, d \xi,
\end{align}
where  $B(t) = \{\xi \in \RR^3: |\xi| \leq g(t) \}$. We now need a pointwise estimate for $||\xi| ^{\frac{1}{2}} \widehat{u} (\xi, t)| ^2$ in $B(t)$, so we take $\Lambda ^{\frac{1}{2}}$ of the solution obtained as a fixed point of \eqref{eqn:fixed-point-operator} and then Fourier transform to have
\begin{align}
\label{eqn:solution-in-ball}
\int _{B(t)} ||\xi| ^{\frac{1}{2}}  \widehat{u} (\xi, t)|^2 \, d \xi  & \leq  C \int _{B(t)} |e^{-t|\xi|^2} |\xi| ^{\frac{1}{2}}  \widehat{u _0} (\xi)|^2 \, d \xi \notag \\ & + C \int _{B(t)} \left( \int_0 ^t e^{-(t-s)|\xi|^2} |\xi| ^{\frac{1}{2}}  {\mathcal F} \left( \mathbb{P}\nabla \cdot  (u \otimes u) \right) (\xi,s) \, ds \right) ^2 \, d \xi.
\end{align}
The estimate for the linear part is
\begin{align}
\label{eqn:decay-linear-ns-h1halfdot}
\int _{B(t)} |e^{-t|\xi|^2} |\xi| ^{\frac{1}{2}}  \widehat{u _0} (\xi)|^2 \, d \xi & \leq   \int _{\RR ^3} |e^{-t|\xi|^2} |\xi| ^{\frac{1}{2}} \widehat{u_0} (\xi)|^2 \, d \xi \leq C (1+t) ^{- \left( \frac{3}{2} + q^{\ast} \right)},
\end{align}
where we used  Theorem \ref{characterization-decay-l2} and $q^{\ast} =  r^{\ast} \left(\Lambda ^{\frac{1}{2}} u_0 \right)$.

We address now the nonlinear part. By the Fubini Theorem,
\begin{align*}
&\int _{B(t)}  \left( \int_0 ^t e^{-(t-s)|\xi|^2}  |\xi| ^{\frac{1}{2}}  {\mathcal F} \left( \mathbb{P}\nabla \cdot  (u \otimes u) \right) (\xi,s) \, ds \right) ^2 \, d \xi \\ &\leq \int_{B(t)} |\xi|\left(\int_0^t  |{\mathcal F} \left( \mathbb{P}\nabla \cdot  (u \otimes u) \right) (\xi,s)| \, ds\right)^2d\xi\\
&=\int_{B(t)}|\xi|\int_0^t\int_0^t|{\mathcal F} \left( \mathbb{P}\nabla \cdot  (u \otimes u) \right) (\xi,s)||{\mathcal F} \left( \mathbb{P}\nabla \cdot  (u \otimes u) \right) (\xi,s')| \, \,dsds'd\xi\\
&\leq g(t)\int_0^t\int_0^t\int_{B(t)}|{\mathcal F} \left( \mathbb{P}\nabla \cdot  (u \otimes u) \right) (\xi,s)| \,|{\mathcal F} \left( \mathbb{P}\nabla \cdot  (u \otimes u) \right) (\xi,s')| \,d\xi dsds'\\
&\leq g(t)\int_0^t\int_0^t\left(\int_{B(t)}|{\mathcal F} \left( \mathbb{P}\nabla \cdot  (u \otimes u) \right) (\xi,s)|^2d\xi\right)^{\frac{1}{2}}\left(\int_{B(t)}|{\mathcal F} \left( \mathbb{P}\nabla \cdot  (u \otimes u) \right) (\xi,s')|^2d\xi\right)^{\frac{1}{2}}dsds'.
\end{align*}
Next, we observe that
\[
|{\mathcal F} \left( \mathbb{P} \nabla \cdot (u \otimes u) \right)|  \leq C |\xi|\, | \mathcal{F}(|u|^2)|,
\]
so, by using H\"older and Hausdorff-Young inequalities
\begin{align*}
    \int_{B(t)}|{\mathcal F} \left( \mathbb{P}\nabla \cdot  (u \otimes u) \right) (\xi,s)|^2d\xi
    & \leq  C   \int_{B(t)}|\xi|^2 | \mathcal{F}(|u|^2)|^2d\xi\\
    & \leq C\left(\int_{B(t)}|\xi|^6 d\xi\right)^{\frac{1}{3}}\left(\int_{B(t)}|\mathcal{F}(|u|^2)|^3\right)^{\frac{2}{3}}\\
    &\leq C g(t)^3 \| |u|^2 \|_{L^{\frac{3}{2}}}^2 = C g(t)^3 \| u \|_{L^3}^4 \\
    & \leq  C g(t)^3  \|u\|_{\dot{H}^{\frac{1}{2}}}^4,
\end{align*}
where we have used that $\dot{H} ^{\frac{1}{2}} (\RR ^3) \subset L^3 (\RR ^3)$.
Plugging this in the previous estimate we arrive at
\begin{align}
\label{eqn:decay-nonlinear-ns-h1halfdot}
&\int _{B(t)}  \left( \int_0 ^t e^{-(t-s)|\xi|^2}  |\xi| ^{\frac{1}{2}}  {\mathcal F} \left( \mathbb{P}\nabla \cdot  (u \otimes u) \right) (\xi,s) \, ds \right) ^2 \, d \xi \nonumber \\
& \leq C g(t)^4   \left( \int _0 ^t \Vert u (s) \Vert _{\dot{H}^{\frac{1}{2}}} ^2  \, ds\right)^2 \nonumber\\
& \leq C g(t)^4   t \int _0 ^t \Vert u (s) \Vert _{\dot{H}^{\frac{1}{2}}} ^4  \, ds.
\end{align}

We first choose $g^2(t) =  \frac{3}{2}  \left( (e+t) \ln (e+t) \right) ^{- 1}$. From \eqref{eqn:key-inequality-ns},  \eqref{eqn:solution-in-ball}, \eqref{eqn:decay-linear-ns-h1halfdot} and \eqref{eqn:decay-nonlinear-ns-h1halfdot}, as we know  that $\Vert u (s) \Vert _{\dot{H}^{\frac{1}{2}}} \leq C$ for all $ s > 0 $, we obtain
\begin{align}
\label{eqn:key-inequality-h1halfdot}
\frac{d}{dt} &  \left( \left( \ln (e+t) \right)^{ 3}  \Vert u (t) \Vert _{\dot{H} ^{\frac{1}{2}}} ^2 \right) \notag \\ & \leq  C \left( (e+t) \ln (e+t) \right) ^{-1}  \left( \ln (e+t) \right)^{ 3} \int _{B(t)} ||\xi| ^{\frac{1}{2}} \widehat{u} (\xi, t)| ^2 \, d \xi \notag \\ & \leq C \frac{\left( \ln (e+t) \right) ^{ 2}}{e+t} \left((1+t) ^{- \left( \frac{3}{2} + q^{\ast} \right)} + \left( \ln (e+t) \right) ^{-2} \right).
\end{align}
We check that
\begin{displaymath}
\int_0 ^t \frac{\left( \ln (e+s) \right) ^{ 2}}{e+s} (1+s) ^{- \left( \frac{3}{2} + q^{\ast} \right)} \, ds \leq C \int_0 ^t \frac{\left( \ln (e+s) \right) ^{ 2}}{(e+s) ^{\left( \frac{5}{2} + q^{\ast} \right)} } \, ds \leq C = C (q^{\ast}).
\end{displaymath}
Thereby, after integrating (\ref{eqn:key-inequality-h1halfdot}), we arrive at the first estimate
\begin{displaymath}
\Vert u (t) \Vert _{\dot{H} ^{\frac{1}{2}}} ^2 \leq C \left( \ln (e+t) \right) ^{ - 2},
\end{displaymath}
which we use to bootstrap in order to obtain finer decay estimates.  Now, let $g^2(t) = \frac{\alpha}{2}(1+t) ^{- 1}$, for $\alpha > 0$ large enough. Using  \eqref{eqn:key-inequality-ns}, \eqref{eqn:decay-linear-ns-h1halfdot} and \eqref{eqn:decay-nonlinear-ns-h1halfdot} again we get
\begin{align}
\label{eqn:first-equation}
\frac{d}{dt} &  \left( \left( 1+t \right)^{\alpha}  \Vert u (t) \Vert _{\dot{H} ^{\frac{1}{2}}} ^2 \right) \notag \\ & \leq C (1+t) ^{\alpha-1}  \int _{B(t)} ||\xi| ^{\frac{1}{2}} \widehat{u} (\xi, t)| ^2 \, d \xi \notag \\ & \leq C (1+t) ^{\alpha - 1} \left( (1+t) ^{- \left( \frac{3}{2} + q^{\ast} \right)} + (1+t) ^{-1} \int_0 ^t   \frac{\Vert u (s) \Vert _{\dot{H} ^{\frac{1}{2}}} ^2}{\left( \ln (e+s) \right) ^{ 2}} \, ds  \right).
\end{align}
Choose $\alpha$ such that  $\alpha > \max \left\{\frac{3}{2} + q^{\ast}, 1 \right\}$. After integrating and dividing by $(1+t)^{\alpha-1}$ on both sides we obtain
\begin{equation}
\label{eqn:second-equation}
( 1+t)  \Vert u (t) \Vert _{\dot{H} ^{\frac{1}{2}}} ^2 \leq  C (1+t) ^{- \left( \frac{1}{2} + q^{\ast} \right)} + C \int_0 ^t   \frac{(1+s) \Vert u (s) \Vert _{\dot{H} ^{\frac{1}{2}}} ^2}{(1+s) \left( \ln (e+s) \right) ^{2}}  ds .
\end{equation}

Suppose $q^{\ast}  \leq - \frac{1}{2}$.  We then apply Gronwall's inequality from Proposition \ref{gronwall-2} with
\begin{align}
\psi(t) & = ( 1+t)  \Vert u (t) \Vert _{\dot{H} ^{\frac{1}{2}}} ^2, \quad a(t) = C (1+t) ^{- \left( \frac{1}{2} + q^{\ast} \right)} \notag \\  &  k (t) = \frac{C}{(1+t) \left( \ln (e+t) \right) ^{ 2}}.  \notag
\end{align}
As
\begin{displaymath}
\int _0 ^t  k(s)  \, ds = C\int _0 ^t \frac{ds}{(1 + s) \left( \ln (e+s) \right) ^{ 2}} \leq C,
\end{displaymath}
we obtain
 \begin{displaymath}
  ( 1+t)  \Vert u (t) \Vert _{\dot{H} ^{\frac{1}{2}}} ^2  \leq C (1+t) ^{- \left( \frac{1}{2} + q^{\ast} \right)}
\end{displaymath}
which leads, to
\begin{equation}
\label{case1}
\Vert u (t) \Vert _{\dot{H} ^{\frac{1}{2}}} ^2  \leq C (1+t) ^{- \left( \frac{3}{2} + q^{\ast} \right)}.
\end{equation}

We now address the case $q^{\ast} > -  \frac{1}{2}$.  We proceed as before,  but instead of obtaining (\ref{eqn:second-equation}) from (\ref{eqn:first-equation}), we get
\begin{displaymath}
\Vert u (t) \Vert _{\dot{H} ^{\frac{1}{2}}} ^2  \leq  C(1+t) ^{- \left( \frac{3}{2} + q^{\ast} \right)} + C(1+t) ^{-1} \int_0 ^t   \frac{\Vert u (s) \Vert _{\dot{H} ^{\frac{1}{2}}} ^2}{ \left( \ln (e+s) \right) ^{ 2}} \, ds.  \notag
\end{displaymath}
We now use Proposition \ref{gronwall-1} with
\begin{align}
x(t) & = \Vert u (t) \Vert _{\dot{H} ^{\frac{1}{2}}} ^2,  \quad a(t) = C (1+t) ^{- \left( \frac{3}{2} + q^{\ast} \right)} \notag \\  & b(t) = \frac{C}{1 + t}, \quad k(t) = \frac{1}{ \left( \ln (e+t) \right) ^{ 2}}.  \notag
\end{align}
As
\begin{displaymath}
\int _s ^t b(r) k(r) \, dr = C\int _s ^t \frac{dr}{(1 + r)  \left( \ln (e+r) \right) ^{ 2}} \leq  C,
\end{displaymath}
hence we obtain
 \[
\Vert u (t) \Vert _{\dot{H} ^{\frac{1}{2}}} ^2  \leq C (1+t) ^{- \left( \frac{3}{2} + q^{\ast} \right)}  + C (1+t) ^{-1} \int_0 ^t \frac{(1+s) ^{- \left( \frac{3}{2} + q^{\ast} \right)}}{ \left( \ln (e+s) \right) ^{ 2}} \, ds.
\]
The condition $q^{\ast} > -  \frac{1}{2}$ implies that the integral on the right hand side is finite,  so  we have
\begin{displaymath}
\Vert u (t) \Vert _{\dot{H} ^{\frac{1}{2}}} ^2   \leq C (1+t) ^{- \left( \frac{3}{2} + q^{\ast} \right)} + C (1+t) ^{-1}.
\end{displaymath}
Thus,
\begin{displaymath}
\Vert u (t) \Vert _{\dot{H} ^{\frac{1}{2}}} ^2  \leq C (1+t) ^{-1}.
\end{displaymath}
This decay together with \eqref{case1}  entails the desired result.  $\Box$

\section{Navier-Stokes-Coriolis equations} \label{NSC}
We first describe  results concerning Navier-Stokes-Coriolis equations \eqref{eqn:navier-stokes-coriolis-omega} which we will need to prove our estimates. We closely follow Hieber and Shibata \cite{MR2609321} (see also Iwabuchi and Takada \cite{MR3096523}).

\subsection{Linear part} \label{linear_part} We first address the linear part of \eqref{eqn:navier-stokes-coriolis-omega}, i.e.
\begin{align}
\label{eqn:parte-lineal}
\partial _t u +  \Omega e_3 \times u + \nabla p & =  \Delta u , \nonumber \\ div \, u & =  0,  \\ u(x,0) & = u_0 (x). \nonumber
\end{align}
Taking the Fourier transform of \eqref{eqn:parte-lineal} we can explicitly solve this equation and see that
\begin{align}
\label{eqn:solution-linear-part}
\widehat{u} (\xi, t) & = e^{t \mathcal{M} (\xi)} \widehat{u}_0 \nonumber \\ & = \cos \left( \Omega \frac{\xi_3}{|\xi|} t \right) e^{- t|\xi|^2 } I_3 \widehat{u}_0 (\xi)  + \sin \left( \Omega \frac{\xi_3}{|\xi|} t \right) e^{-t |\xi|^2 } R(\xi) \widehat{u}_0 (\xi),
\end{align}
and, as a consequence,
\begin{displaymath}
\widehat{p} (\xi) = i \frac{\Omega}{|\xi|^2} \left( \xi_2 \widehat{u_1} (\xi) - \xi_1 \widehat{u_2} (\xi) \right).
\end{displaymath}
In \eqref{eqn:solution-linear-part}, $I_3$ is the identity matrix and
\begin{displaymath}
R (\xi) =  \left( \begin{array}{ccc} 0 & \frac{\xi _3}{|\xi|} &  - \frac{\xi_2}{|\xi|} \\   - \frac{\xi_3}{|\xi|} & 0  &  \frac{\xi _1}{|\xi|} \\  \frac{\xi_2}{|\xi|} &  - \frac{\xi_1}{|\xi|} & 0 \end{array} \right).
\end{displaymath}
We then obtain a semigroup such that for divergence free $f \in  L^2  (\RR ^3) $
\begin{equation}
\label{eqn:semigroup-coriolis}
T(t) \, f = \mathcal{F} ^{-1} \left( \cos \left( \Omega \frac{\xi_3}{|\xi|} t \right) e^{- t|\xi|^2 } I_3 \widehat{f} (\xi) + \sin \left( \Omega \frac{\xi_3}{|\xi|} t \right) e^{-t |\xi|^2 } R(\xi) \widehat{f} (\xi) \right).
\end{equation}
In fact, by Mikhlin's Theorem, $T(t)$ can be extended to a $C_0$ semigroup $T_p(t)$ on divergence free $ L^p  (\RR ^3) $, where $1 < p < \infty$.

We  describe now some pointwise properties of $\widehat{u} (\xi,t)$ we need for using the Fourier Splitting Method.  Through an elementary computation there follows
\begin{displaymath}
|\widehat{u} (\xi, t)| = |e^{t \mathcal{M} (\xi)} \widehat{u}_0| \leq C \, e^{-t |\xi|^2} |\widehat{u}_0|,
\end{displaymath}
from which, using Theorem \ref{characterization-decay-l2}, we obtain that for $u_0 \in L^2(\RR^3)$ with decay character $- \frac{3}{2} <  r^{\ast} =  r^{\ast} (u_0) < \infty$
\begin{equation*}
%\label{eqn:l2-decay-linear-coriolis-operator}
\Vert u(t) \Vert _{L^2} ^2  = \Vert e^{t \mathcal{M} (\xi)} \widehat{u}_0 \Vert _{L^2} ^2  \leq C (1 + t) ^{- \left( \frac{3}{2} + r^{\ast}  \right)}.
\end{equation*}

\subsection{Existence of solutions in $\dot{H}^{\frac{1}{2}} (\RR ^3)$} \label{section-existence-hdot} To obtain existence of global solutions to \eqref{eqn:navier-stokes-coriolis-omega} in the critical, homogeneous space $\dot{H}^{\frac{1}{2}} (\RR ^3)$  we turn \eqref{eqn:navier-stokes-coriolis-omega} into an integral equation
\begin{equation} \label{eqn:integral-equation}
\mathcal{T} (u) = T(t) u_0 - \int_0 ^t T(t-s) \, \mathbb{P} \nabla \cdot\left( u \otimes u  \right)(s) \, ds,
\end{equation}
where $T $ is given by \eqref{eqn:semigroup-coriolis} and find an appropriate Banach space $X$ on which this has a fixed point. We note that there are differences with the Navier-Stokes'  case that must be taken into account.
As we mentioned in Section \ref{NS}, there exists $\epsilon > 0$, such that for initial data with $\Vert u_0 \Vert _{\dot{H} ^\frac{1}{2}} < \epsilon$, the Navier-Stokes equations \eqref{eqn:NS} have a global solution in $C \bigl([0, \infty); \dot{H}^{\frac{1}{2}} (\RR ^3) \bigr) $.  For the Navier-Stokes-Coriolis equations  \eqref{eqn:navier-stokes-coriolis-omega}, Iwabuchi and Takada \cite{MR3096523} and Koh, Lee and Takada \cite{MR3229600} proved that  in $\dot{H} ^s (\RR ^3)$, with $\frac{1}{2} < s < \frac{9}{10}$, the condition
\begin{displaymath}
\Vert u_0 \Vert _{\dot{H}^s} \leq C |\Omega| ^{\frac{1}{2} \left(s - \frac{1}{2} \right)}
\end{displaymath}
implies existence of solutions in $L^{\infty} \bigl([0, \infty);  \dot{H}^s (\RR ^3) \bigr)$. Hence, in both cases,  conditions on the size of initial data lead to global existence.
However, this is  not the case in $\dot{H}^{\frac{1}{2}} (\RR ^3)$, where the conditions for global existence depend on the initial data, not merely its size, as can be seen in the following Theorem.

\begin{Theorem}[Theorem 1.3, Iwabuchi and Takada \cite{MR3096523}] %\label{thm-first-iwabuchi-takada}
For any $u_0 \in \dot{H}^{\frac{1}{2}} (\RR ^3)$ with $\div \, u_0 = 0$, there exists a $\omega = \omega (u_0)$ such that for $\omega< |\Omega|$ there exists a unique global solution to  \eqref{eqn:navier-stokes-coriolis-omega} in $ C \bigl([0, \infty); \dot{H}^{\frac{1}{2}} (\RR ^3) \bigr) \cap L^4 \bigl([0, \infty); \dot{H}^{\frac{1}{2}} _3 (\RR ^3) \bigr) $.
\end{Theorem}
 Here, $\dot{H} ^s _p = \{ f \in \mathcal{S} ' \left( \RR ^n\right)/P(\RR^n): \Lambda ^s f \in L^p(\RR^n)\}$ is the homogeneous potential space.

\begin{Remark}%\label{profileremark}
The fact that the size condition for global existence involves the profile of initial data instead of of its norm is also usually seen in dispersive critical equation, such as in the $L^2$-critical Schr\"odinger equation, see Cazenave and Weissler \cite{MR1021011} and the generalized critical KdV equation,  see Birnir, Kenig, Ponce, Svanstedt and Vega \cite{MR1396718} (see also Theorem 5.3 and Theorem 7.7 in Linares and Ponce \cite{MR2492151}).
\end{Remark}

In order to be able to use the Fourier Splitting Method to prove decay of solutions we need to show that for a given constant $C > 0$ we have that
\begin{equation}
\label{eqn:estimate-hdotnorm-c}
\Vert u(t) \Vert _{\dot{H} ^{\frac{1}{2}}} < \frac{1}{C}, \qquad \forall \, t > 0.
\end{equation}
We claim that this is possible for $u_0$ with small enough $\dot{H} ^{\frac{1}{2}}$-norm.

We now describe the setting in which Iwabuchi and Takada \cite{MR3096523} prove existence of solutions and we then show how we can attain \eqref{eqn:estimate-hdotnorm-c}.

Consider the ball
\begin{displaymath}
B = \Bigl\{u \in C^0 \bigl([0, \infty), \dot{H} ^{\frac{1}{2}}  \bigr), \Vert u \Vert _{L^4 _t \dot{H} ^{\frac{1}{2}} _{3,x}} \leq 2 \, \delta, div \, u = 0  \Bigr\}
\end{displaymath}
in the metric space $X = L^4 \bigl([0, \infty);  \dot{H} ^{\frac{1}{2}} _3 (\RR^3)\bigr)$ with the metric induced by the norm and where $ \delta > 0 $ is  small enough. The integral operator  associated to equation \eqref{eqn:navier-stokes-coriolis-omega} is a contraction in $B$, and hence has a fixed point. The fixed point argument does not furnish the smallness of the solution in $L^\infty\bigl((0,\infty);\dot{H} ^{\frac{1}{2}}(\RR^3) \bigr)$. However,
the solution obtained also obeys the key estimate
\begin{equation}
\label{eqn:estimate-small}
\sup _{t > 0} \Vert u(t)  \Vert _{\dot{H} ^{\frac{1}{2}}} \leq   C'  \Vert u_0 \Vert _{\dot{H} ^{\frac{1}{2}}} + 4 C ''   \delta ^2,
\end{equation}
where $C', C'' > 0$ do not depend on $\delta$, see (3.5), page 739 in  \cite{MR3096523}. We then see that \eqref{eqn:estimate-hdotnorm-c} holds provided $\Vert u_0 \Vert _{\dot{H} ^{\frac{1}{2}}}$ and $\delta > 0$ are small enough.

\subsection{Norm in $\dot{H}^{\frac{1}{2}} (\RR ^3)$ is a Lyapunov function} \label{sub-sub-lyapunov} In order to prove that the $\dot{H}^{\frac{1}{2}} (\RR ^3)$ is a Lyapunov function also for the Navier-Stokes-Coriolis  equations, we need the following Lemma, which we state with no proof.

\begin{Lemma}
Let $w \in C\left([0,T]; \mathcal{S} '(\RR ^n) \right)$ be the solution to
\begin{align*}
\partial _t w + L w - \Delta w  & =   f,  \nonumber \\ w(x,0) & = w_0 (x) ,
\end{align*}
where $f \in L^2 \left([0,T]; \dot{H} ^{s-1}(\RR^n) \right)$, $w_0 \in \dot{H}^{s} (\RR^n)$ and $L:\mathcal{S} ' (\RR ^n) \to \mathcal{S} ' (\RR ^n)$ is a continuous linear map induced by a linear map $\widetilde{L}:\mathcal{S} (\RR ^n) \to \mathcal{S} (\RR ^n)$. If $\langle Lw, w \rangle _{\dot{H} ^s} = 0$, then we have the energy identity
\begin{equation*}
%\label{eqn:energy-identity}
\Vert w(t) \Vert ^2 _{\dot{H} ^s}  + 2 \int _0 ^t \Vert \nabla w(\tau) \Vert ^2  _{\dot{H}^s} \, d\tau = \Vert w_0 \Vert ^2 _{\dot{H} ^s} + 2  \int _0 ^t \langle w (\tau), f(\tau) \rangle _{\dot{H}^s} \, d \tau.
\end{equation*}
\end{Lemma}
We now project equation \eqref{eqn:navier-stokes-coriolis-omega} onto the space of divergence free vector fields with $\mathbb{P}$. If we consider $\widetilde{L} (u) = \mathbb{P} \left( e_3 \times u \right)$, as a consequence of
\begin{displaymath}
\widehat{\mathbb{P} \left( e_3 \times u \right)} = \left( - \left(1 - \frac{\xi_1 ^2}{|\xi|^2} \right) \widehat{u_2} - \frac{\xi_1 \xi_2}{|\xi|^2} \widehat{u_1}, \frac{\xi_1 \xi_2}{|\xi|^2} \widehat{u_2} + \left(1 - \frac{\xi_2 ^2}{|\xi|^2} \right) \widehat{u_1},  \frac{\xi_1 \xi_3}{|\xi|^2} \widehat{u_2} -  \frac{\xi_3 \xi_2}{|\xi|^2} \widehat{u_1} \right),
\end{displaymath}
and $\widehat{div \, u} = \xi \cdot \widehat{u} = 0$,  we have that $\langle L u , u  \rangle _{\dot{H}^s} = \langle \mathbb{P} \left( e_3 \times u \right), u  \rangle _{\dot{H}^s} = 0$. We then, as before, obtain the energy identity
\begin{align}
\Vert u(t) \Vert ^2 _{\dot{H} ^{\frac{1}{2}}}   + 2 \int _0 ^t \Vert \nabla u(\tau) \Vert ^2  _{\dot{H}^{\frac{1}{2}}} \, d\tau &  = \Vert u_0 \Vert ^2 _{\dot{H} ^{\frac{1}{2}}} \notag \\ & - 2  \int _0 ^t \langle  u(\tau), \mathbb{P} \left( u \cdot \nabla u (\tau) \right) \rangle _{\dot{H}^{\frac{1}{2}}} \, d \tau, \notag
\end{align}
and through Proposition 5.13 in Bahouri {\em et al.}  \cite{MR2768550}, we get that the $\dot{H} ^{\frac{1}{2}}$-norm is a Lyapunov function for the Navier-Stokes-Coriolis equations.

\subsection{Proof of Theorem \ref{thm-decay-h-dot-one-half}.} From the previous Section, we know that the $\dot{H} ^{\frac{1}{2}}$-norm is a Lyapunov function, so proceeding as in the proof of Theorem \ref{decay-ns-critical} we obtain
\begin{displaymath}
\frac{1}{2} \frac{d}{dt} \Vert u(t) \Vert _{\dot{H} ^{\frac{1}{2}}} = - \left(1 -  C \Vert  u(t) \Vert _{\dot{H} ^{\frac{1}{2}}} \right) \Vert \nabla u(t) \Vert _{\dot{H} ^{\frac{1}{2}}}.
\end{displaymath}
From Section \ref{section-existence-hdot}, by choosing  $\Vert u_0 \Vert _{\dot{H}^{\frac{1}{2}}}$ and $\delta > 0$ small enough,  \eqref{eqn:estimate-small} leads to $\Vert u(t) \Vert _{\dot{H} ^{\frac{1}{2}}} < \frac{1}{C}$, hence we can use the Fourier Splitting Method.  Following the proof of Theorem \ref{decay-ns-critical} we arrive at the key estimate  \eqref{eqn:key-inequality-ns}.

Hence, we need a pointwise estimate for
\begin{align} \label{eqn:pointwise}
\int _{B(t)} & ||\xi| ^{\frac{1}{2}}  \widehat{u} (\xi, t)|^2 \, d \xi \notag \\ & \leq  C \int _{B(t)} |e^{t \mathcal{M} (\xi)} |\xi| ^{\frac{1}{2}}  \widehat{u _0}|^2 \, d \xi  \notag\\ & + C \int _{B(t)} \left( \int_0 ^t e^{(t -s) \mathcal{M} (\xi)} |\xi| ^{\frac{1}{2}}  {\mathcal F} ( \mathbb{P} \nabla \cdot \left( u \otimes u \right)) (\xi,s) \, ds \right) ^2 \, d \xi.
\end{align}
The estimate for the linear part is
\begin{align*}
%\label{eqn:decay-linear-h1halfdot}
\int _{B(t)} |e^{t \mathcal{M}(\xi)} |\xi| ^{\frac{1}{2}} \widehat{u_0}|^2 \, d \xi & \leq \int _{\RR ^3} |e^{t \mathcal{M}(\xi)} |\xi| ^{\frac{1}{2}} \widehat{u_0}|^2 \, d \xi \notag \\ & \leq C \int _{\RR ^3} |e^{-t|\xi| ^2} |\xi| ^{\frac{1}{2}} \widehat{u_0}|^2 \, d \xi \leq C (1+t) ^{- \left( \frac{3}{2} + q^{\ast} \right)},
\end{align*}
where we used  Theorem \ref{characterization-decay-l2} and $q^{\ast} =  r^{\ast} \left(\Lambda ^{\frac{1}{2}} u_0 \right)$.

For the nonlinear part we proceed as in \eqref{eqn:decay-nonlinear-ns-h1halfdot}
to conclude
\begin{align*}
& \int_{B(t)}\left(  \int _0 ^t e^{(t -s)  \mathcal{M} (\xi)} |\xi| ^{\frac{1}{2}}  {\mathcal F} ( \mathbb{P} \nabla \cdot \left( u \otimes u \right)) (\xi,s) \, ds \right) ^2 \\
&\leq \int _{B(t)}  \left( \int_0 ^t e^{-(t-s)|\xi|^2}  |\xi| ^{\frac{1}{2}}  |{\mathcal F} \left( \mathbb{P}\nabla \cdot  (u \otimes u) \right) (\xi,s)| \, ds \right) ^2 \, d \xi \\
& \leq C g(t)^4   t \int _0 ^t \Vert u (s) \Vert _{\dot{H}^{\frac{1}{2}}} ^4  \, ds.
\end{align*}

 We can now continue in the same way as in the proof of Theorem \ref{decay-ns-critical} to obtain the corresponding estimate.  $\Box$

 \subsection{Additional decay estimate in $H^{\frac{1}{2}} (\RR ^3)$}
We now turn our attention to the space $H ^{\frac{1}{2}} (\RR ^3)$, where existence of solutions was proved by Hieber and Shibata \cite{MR2609321}  by adapting arguments of Fujita and Kato \cite{MR166499}.   They managed to express the conditions for existence of solutions in terms of the size of $\Vert u_0 \Vert _{H ^{\frac{1}{2}}}$ only, with no reference to the rotation strength $\Omega$. Indeed,
a solution to  \eqref{eqn:navier-stokes-coriolis-omega} will be a fixed point of the map $\mathcal{T}:E \to E$ given by \eqref{eqn:integral-equation}. In this case, we have $E = E_{\epsilon}$, where
\begin{displaymath}
E_{\epsilon}  =  \left.
\begin{cases}
u \in C \bigl(  [0, \infty); H ^{\frac{1}{2}}  (\RR ^3)  \bigr) :  <u>_t \, \leq \, \epsilon, \forall t > 0;  \, \lim _{t \to 0^+} \Vert u(t) - u_0 \Vert _{H^{\frac{1}{2}}} = 0, \\
u \in C \bigl(  (0, \infty); L^q  (\RR ^3)  \bigr):  \lim _{t \to 0^+} [u]_{q, \frac{1}{2}- \frac{3}{2q}, t} = 0, \, div \, u =0 \\
\nabla u \in C \bigl(  (0, \infty); L^2  (\RR ^3)  \bigr):  \lim _{t \to 0^+} [\nabla u]_{2, \frac{1}{4}, t} = 0,
\end{cases}
\right\}
\end{displaymath}
with metric
\begin{align*}
<u>_t & = \sup _{0 < s < t} \Vert u(s) \Vert _{H^{\frac{1}{2}}} + \sup _{0 < s < t} s ^{\frac{1}{2} - \frac{3}{2q}} \Vert u(s) \Vert _{L^q} \notag \\ & + \sup _{0 < s < t} s ^{\frac{3}{4} - \frac{3}{2q}} \Vert u(s) \Vert _{L^q} + \sup _{0 < s < t} s^{\frac{1}{4}} \Vert \nabla u(s) \Vert _{L^2}.
\end{align*}
The map $\mathcal{T}$ has a fixed point in $E_{\epsilon}$, provided  $u_0 \in H ^{\frac{1}{2}}(\RR^3)$ with $ div \,u_0 =0 $  and $\epsilon > 0$ is small enough, see pages 488 - 490 in Hieber and Shibata \cite{MR2609321}.  This guarantees existence of global solutions in $H^{\frac{1}{2}}(\RR^3)$.
Note that, unlike what happens when studying existence of weak solutions (see Chemin \emph{et al.}  \cite{MR1935994}, \cite{MR2228849}), $\epsilon > 0$ does not depend on $\Omega$.
However, their solutions are just bounded,  i.e. $\Vert u(t) \Vert _{H ^{\frac{1}{2}}} \leq C$,  and no decay is provided.  Our next result establishes such rates.

\begin{Theorem} %\label{thm-decay-h-one-half}
Let $u_0 \in H ^{\frac{1}{2}} (\RR ^3)$, with $ div \, u_0 = 0 $ and $\Vert u_0  \Vert _{ H ^{\frac{1}{2}} } < \epsilon$, for small enough $\epsilon >  0$ and $-\frac{3}{2} < r^{\ast} = r^{\ast} (u_0) < \infty$. Then, for any mild solution to \eqref{eqn:navier-stokes-coriolis-omega} we have that

\begin{displaymath}
\Vert u (t) \Vert _{L^2} ^2 \leq C (1 + t) ^{- \min \{ \frac{3}{2} + r^{\ast}, \frac{5}{2}\}}.
\end{displaymath}
Also,

\begin{displaymath}
\Vert u (t) \Vert _{\dot{H} ^{\frac{1}{2}}} ^2 \leq C (1 + t) ^{- \left( \frac{1}{2} + \min \{ \frac{3}{2} + r^{\ast}, \frac{5}{2}\} \right)}.
\end{displaymath}
\end{Theorem}
Observe that, as is expected,  the  $H ^{\frac{1}{2}}$ decay of solutions with initial datum in $ H ^{\frac{1}{2}}(\RR^3) $ is driven by the $L^2$ decay. Moreover, there is a gain of $\frac{1}{2}$ on the decay in analogy with that happens  for  Navier-Stokes equations.

\begin{proof} First, we observe that  solutions to \eqref{eqn:navier-stokes-coriolis-omega}  satisfy the energy inequality
\begin{displaymath}
\frac{1}{2} \frac{d}{dt} \Vert u(t) \Vert _{L^2} ^2 \leq - \Vert \nabla u(t) \Vert _{L^2} ^2.
\end{displaymath}
As pointed out in Section \ref{linear_part}, we have in frequency space the pointwise estimate
\begin{displaymath}
 |e^{t \mathcal{M} (\xi)} \widehat{f}| \leq C \, e^{-t |\xi|^2} |\widehat{f}|,
\end{displaymath}
so the proof of decay for the $L^2$-norm is analogous as that for the Navier-Stokes equations.

We will use the decay for the $L^2$-norm to prove the decay of solutions in the $\dot{H} ^{\frac{1}{2}}$-norm.
As the $\dot{H} ^{\frac{1}{2}}$-norm is a nonincreasing function (Section  \ref{sub-sub-lyapunov}), proceeding as in the proof of Theorem \ref{thm-decay-h-dot-one-half},  choosing $\epsilon$ small enough we then have
\begin{displaymath}
\frac{d}{dt} \Vert u(t) \Vert^2 _{\dot{H} ^{\frac{1}{2}}} \leq  - 2  \Vert \nabla u(t) \Vert^2 _{\dot{H} ^{\frac{1}{2}}},
\end{displaymath}
thus we can use the Fourier Splitting Method from Section \ref{fourier-splitting-section}.
 Following the proof of Theorem \ref{decay-ns-critical} we arrive at the key estimate  \eqref{eqn:key-inequality-ns}.
 So, we need to  estimate  $||\xi| ^{\frac{1}{2}} \widehat{u} (\xi, t)| ^2$ in $B(t)$ from \eqref{eqn:pointwise}.

For the linear part in \eqref{eqn:pointwise}, it holds that
\begin{align}
\label{eqn:decay-linear}
\int _{B(t)} |e^{t\mathcal{M}(\xi) } |\xi| ^{\frac{1}{2}} \widehat{u_0}|^2 \, d \xi &   \leq \int _{\RR ^3} |e^{t\mathcal{M}(\xi) } |\xi| ^{\frac{1}{2}} \widehat{u_0}|^2 \, d \xi  \notag \\& \leq C\int _{\RR ^3} |e^{-t|\xi|^2 } |\xi| ^{\frac{1}{2}} \widehat{u_0}|^2 \, d \xi   \leq C (1+t) ^{- \left( \frac{1}{2} + \frac{3}{2} + r^{\ast} \right)},
\end{align}
where we used Theorem \ref{characterization-decay-hs}.

For the nonlinear part in \eqref{eqn:pointwise}, by estimating
\begin{equation*}
\label{eqn:nonlinear-term}
|{\mathcal F}( \mathbb{P} \nabla \cdot\left( u \otimes u \right))| \leq C |\mathcal{F} \left( \nabla \cdot (u \otimes u ) \right)| \leq  C |\xi|  |\widehat{u \otimes u}| \leq C |\xi| \Vert u(t) \Vert _{L^2} ^2,
\end{equation*}
there follows
\begin{align*}
\label{eqn:inequality-integral-term}
\left( \int _0 ^t e^{(t - s)|\mathcal{M}(\xi)} |\xi| ^{\frac{1}{2}} {\mathcal F}( \mathbb{P} \nabla \cdot\left( u \otimes u \right)) (\xi,s) \, ds \right)^2 & \leq  C\left( \int _0 ^t e^{-  (t - s) |\xi|^2} \, |\xi| ^{\frac{3}{2}} \Vert u(s) \Vert _{L^2} ^2 \, ds \right)^2 \nonumber \\ & \leq  C |\xi|^{3} \left( \int _0 ^t \Vert u(s) \Vert _{L^2} ^2 \, ds \right)^2.
\end{align*}
 First,  we consider $r^* \neq -\frac{1}{2}$. Using the decay estimate for the $L^2$-norm of $ u $ we obtain
\[\left( \int _0 ^t \Vert u(s) \Vert _{L^2} ^2 \, ds \right)^2 \leq \left \{ \begin{array}{ll}
C , & \text{ if } r^* > -\frac{1}{2}\\
C(1+t)^{-(1+2r^*)} , & \text{ if } r^* < -\frac{1}{2}.
\end{array} \right.
\]
Taking $g^2 (t) = \frac{\alpha}{2 } (1 +t) ^{-1}$ with large enough $\alpha > 0$ leads to
\begin{equation}
\label{eqn:decay_nolinear} \int_{B(t)}|\xi|^3 \left( \int _0 ^t \Vert u(s) \Vert _{L^2} ^2 \, ds \right)^2 d\xi \leq \left \{ \begin{array}{ll}
C(1+t)^{-3} , & \text{ if } r^* > -\frac{1}{2}\\
C(1+t)^{-(4+2r^*)} , & \text{ if } r^* < -\frac{1}{2}.
\end{array} \right.
\end{equation}
If $r^* = -\frac{1}{2}$, we proceed in a slightly different way,  using that
\begin{displaymath}
\left( \int _0 ^t \Vert u(s) \Vert _{L^2} ^2 \, ds \right)^2 \leq C t \int _0 ^t \Vert u(s) \Vert _{L^2} ^4 \, ds  \leq Ct.
\end{displaymath}
Hence
\begin{equation}
\label{eqn:decay_nonlinear_special_case} \int_{B(t)}|\xi|^3 \left( \int _0 ^t \Vert u(s) \Vert _{L^2} ^2 \, ds \right)^2 d\xi \leq C t \int_{B(t)}|\xi|^3 \, d \xi \leq C (1 + t) ^{-2}.
\end{equation}
Now, comparing \eqref{eqn:decay-linear} to \eqref{eqn:decay_nolinear} and \eqref{eqn:decay_nonlinear_special_case}  we notice that the slower decay is given by the nonlinear part, i.e.  \eqref{eqn:decay_nolinear} and \eqref{eqn:decay_nonlinear_special_case},  when $ r^* \geq 1 $ and by the linear part, i.e. \eqref{eqn:decay-linear},  when $r^*<1$. Plugging this in \eqref{eqn:key-inequality-ns}, we then have
\begin{displaymath}
\frac{d}{dt} \left( (1+t) ^{\alpha}  \Vert u (t) \Vert _{\dot{H} ^{\frac{1}{2}}} ^2 \right) \leq C (1+t) ^{\alpha- 1} (1+t) ^{- \left( \frac{1}{2} +  \min\{\frac{3}{2} + r^{\ast}, \frac{5}{2}\} \right)}.
\end{displaymath}
Choosing $\alpha - 1 -  \left( \frac{1}{2} +\min\{\frac{3}{2} + r^{\ast}, \frac{5}{2}\} \right) > - 1$, after integrating, there follows  the desired estimate.
\end{proof}

\end{document}